\documentclass[12 pt]{article}
\newtheorem{theorem}{Theorem}
\newtheorem{definition}{Definition}
\begin{document}
\title{Variational Synthesis\\ of Controlled Dynamic Mappings}
\author{V. Yu. Tertychny-Dauri\\
Saint-Petersburg State University of\\ Information Technologies, Mechanics, and Optics\\
Department of Physics and Engineering\\
49, Kronverkskyi Pr., St-Petersburg. 197101, Russia\\
tertychny-dauri@mail.ru}
\date{}
\maketitle
\begin{abstract}
The article deals with the subject of solving the problem of canonical-map synthesis for
Hamiltonian systems. For this purpose, the controlling-function method has
been developed that allows appropriate changes of the variables
in terms of calculus of canonical variation, starting from their target conditions.
To use the canonical formalism, the initial dynamic system that employs changing
Lagrange multipliers is reduced to a Hamiltonian system in an
expanded phase space, followed by the construction of controlling function.
The algorithm suggested for the canonization of controlled mappings has an
advantage over the known procedures, and first of all, redundancy in the
procedure  that chooses regulated coordinate transformations as a base for
a goal-seeking synthesis scheme.
\end{abstract}
\newpage
\section*{Introduction}
Below we suggest a new procedure of variational changes in dynamical and
mechanical systems, which is named the method of controlled (or synthesized)
mappings. In distinction to the known methods of generative functions, those of
Lie generators, parametrization in the Hilbert-Courant form, etc. [1--5], the method allows
a large extent of generality and versatility (for example, with respect to
the transformations of coordinate axes, etc.) because it makes possible to
synthesize mappings in dependence on the chosen target conditions without
resort to some laborious-intensive procedures (as it takes place in
a~classical variant of the generative-function method).

The formalism of generative functions unfortunately does not permit beforehand
a new system of variables suitable to finding
solution of the selected problems (normalization, integration, and so
on). On the contrary, the method of controlled mapping makes
possible  to solve
these problems constructively basing on the
beforehand formulated conditions. It might be caused by some
shortcomings in the method itself, which presented, as so often was the case, little more than
derivation of the related equations in the required quantities, but did not present their
analytical solution.

In Section 1, a controlled change of variables  in dynamic systems has been
considered. For this purpose, canonical equations are derived
variationally, with giving proof of a~necessary Weierstrass condition for a
minimal action functional in terms of invariant Hilbert integral. The
controlled mapping is synthesized in phase space by means of a~controlling
function, and its discovery becomes the most important problem of the
method.

In Section 2 the most distinctive singularities of the suggested controlled
transformations are pointed. Here we, first of all, reveal the conditions for appearing
the properties of~canonicity and invariance. Various forms of canonicity
criteria for the controlled mappings are found, and some possible
ways of their support are discussed. Great attention is paid to the derivation of
Hamilton-Jacobi equation for the controlling potential. For~the introduced
controlled mappings, the invariance properties are established.

Section 3 presents description of energy transformations in an initial
dynamic system, which were induced by the canonical controlled mappings. The
relationship between the acting and controlling functions is established, which provides the controlling
transformation field in phase space. By a basic target
condition imposed on the new Hamiltonian, they present an algorithmic scheme for
finding the controlling function.

Section 4 is small in volume but important in theoretical generalizations.
The~direct analogy between small controlled canonical mappings and
infinitesimal canonical transformations is firstly revealed. It is secondly  shown that
the canonization of controlled transformations provides target conditions
without using a special choice of initial data.

\section{Controlled change of variables}
The controlled transformation of variables in any moving dynamic system is
presumed to~be a change of system variables in the process of system movement as
time passes, i.~e. supplying a~controlled signal in the form of controlled mapping
to an input converter of initial variables $x(t)$.

Let movement of the dynamic system be described by an $n$-dimensional vector
equation
$$\dot{x} = f(x,t),\qquad x \in R^n,\eqno(1.1)$$
where $x = (x_1,\,...,\,x_n)$  is the phase vector,
$x_i = x_i(t),\,i = \overline{1,n}$ are phase variables (coordinates,
velocities), $t \in [\,t_0,\,t_1\,] \subset
R$ is time. It is thought that $f,\,\partial f/
\partial x,\,\partial f / \partial t \in C^1[\,t_0,\,t_1\,]$  on the right  side of equation (1.1)
satisfy the existence and uniqueness condition applied to the Cauchy
problem in some limited region $D \subset R^n \times R$ , i.~e. equation (1.1) is satisfied by a
unique continuous integral path in the expanded phase space $R^n \times R$ that originates at
the point $x_{0i} = x_i(t_0)$ and terminates at the point $x_{1i} = x_i(t_1),\,i = \overline{1,n}.$

Given differential constraints in the form of equation of motion (1.1), let us
introduce the Lagrange multipliers $\lambda_i(t),\,i = \overline{1,n}, $
where $\lambda_i(t)$ is the unknown temporal functions
determined as part of the solution of optimization problem for the action
functional~$S$, with the provision of constrained equation (1.1). The Lagrange
multipliers can be seen as weighting action-functional integral multipliers
$$S = \int_{t_0}^{t_1} L(x,\dot{x},\lambda,t)\,dt,\eqno(1.2)$$
where $L = \lambda\,(\dot{x} - f)$ is the Lagrange function (Lagrangian).

The variation of variable quantity is, as we know, the difference of values of this
quantity on a permissible reference trajectory and on an appropriate
permissible trajectory called the comparison trajectory.
The variation of the quantity calculated for the same moment of time $t$ is
called the isochronous variation and denoted by $\delta.$ If in calculating the
variation the value of the quantity is considered on the comparison trajectory
at time $t + \Delta t,$ this variation is called a total variation. Let the total
variation be denoted by $\Delta.$

Supply also the quantities taken on the comparison trajectory by a line over.
Then the relationship between the total and isochronous variations is specified by the equalities
$$\Delta x_i = \bar{x}_i(t + \Delta t) - x_i(t) = \bar{x}_i(t + \Delta
t) - \bar{x}_i(t) + \delta x_i,$$
$$\Delta x_i = \delta x_i + \dot{\bar{x}}_i\,\Delta t,\qquad i = \overline{1,n},
\eqno (1.3)$$
and the quantity $\bar{x}_i$ in formula (1.3) substitutes for $\dot{\bar{x}_i}$
if they were considered to be
diverged infinitesimally. It is here thought that we have the permutation
relations $\delta \dot{x}_i = \dot{\bar{x}}_i - \dot{x}_i = d\,(\bar{x}_i - x_i)/dt = d\,
(\delta x_i)/dt.$

To seek $\Delta S$, let us go from the Lagrange function $L$ to the corresponding Hamilton
function $H$ using the Legendre transformation
$$H(x,\lambda,t) = {{\partial L}\over{\partial \dot{x}}}\,\dot{x} - L,\qquad L =
L(x,\dot{x},\lambda,t).\eqno(1.4)$$
According to the definition of function $L$
and by (1.4), we obtain
$${{\partial L}\over{\partial\dot{x}}} = \lambda,\qquad H = \lambda f.\eqno(1.5)$$
The vector multiplication is here and further on considered as a scalar
product.

With relations (1.4), (1.5) we can write [6--8]
$$\Delta S = \int_{t_0+\Delta t_0}^{t_1+\Delta t_1}\,\bar{\lambda}\,\dot{\bar{x}}
\,dt - \int_{t_0}^{t_1}\,\lambda\dot{x}\,dt - \Delta\,\int_{t_0}^{t_1}\,H\,dt\ $$
$$=\ \int_{t_0}^{t_1}\,\bigl(\bar{\lambda}\,\dot{\bar{x}} - \lambda\dot{x}\bigr)\,
dt + \bar{\lambda}\,\dot{\bar{x}}\,\big|_{t_1}\,\Delta t_1 - \bar{\lambda}\,
\dot{\bar{x}}\,\big|_{t_0}\,\Delta t_0\ $$
$$-\ \int_{t_0}^{t_1}\,\delta H\,dt - H_1\,\Delta t_1 + H_0\,\Delta t_0,
\eqno(1.6)$$
where $H_0 = H\,\big|_{t=t_0},\,H_1 = H\,\big|_{t=t_1},\,\bar{\lambda} - \lambda =
\delta\lambda$, and all equalities in (1.6) are valid accurate within infinitesimals
of the order of smallness higher than first. The first summand in the
right side of formula (1.6) may be written with the mentioned accuracy in the
form
$$\int_{t_0}^{t_1}\bigl(\bar{\lambda}\,\dot{\bar{x}} - \lambda\dot{x}\bigr)\,dt =
\int_{t_0}^{t_1} \lambda\,\bigl(\dot{\bar{x}} - \dot{x}\bigr)\,dt + \int_{t_0}^{t_1}
\,\dot{\bar{x}}\,\delta\lambda\,dt\ $$
$$=\ \int_{t_0}^{t_1}\lambda\,{d\over{dt}}\,\bigl(\delta x\bigr)\,dt +
\int_{t_0}^{t_1}\,\dot{x}\,\delta\lambda\,dt = \int_{t_0}^{t_1}\bigl(-\,\dot{\lambda}
\,\delta x + \dot{x}\,\delta\lambda\bigr)\,dt + \lambda\,\delta x\,\big|_{t_0}^
{t_1}.\eqno(1.7)$$

Using relations (1.3), (1.7), and
$$\delta H(x,\lambda,t) = {{\partial H}\over{\partial x}}\,\delta x + {{\partial
H}\over{\partial \lambda}}\,\delta\lambda$$
for the calculation of $\Delta S$ being accurate with the second and higher orders of smallness, we obtain
$$\Delta S = -\,\int_{t_0}^{t_1}\biggl[\,\biggl({{\partial H}\over{\partial x}}
+ \dot{\lambda}\biggr)\,\delta x + \biggl({{\partial H}\over{\partial\lambda}} -
\dot{x}\biggr)\,\delta\lambda\,\biggr]\,dt\ $$
$$+\ \bigl(\lambda\,\Delta x - H\,\Delta t\bigr)\,\big|_{t_0}^{t_1}.\eqno(1.8)$$

If there is an extremal of the functional $S$, it is necessary to fulfil the
stationary condition $\Delta S = 0.$ When the variations in the integral (1.8) are independent
and when we take account to the equalities (1.5), we obtain
necessary conditions in the form of:
\begin{itemize}
\item[(1)] Euler equations (Euler-Lagrange) with respect to the multipliers $\lambda$
$$\dot{\lambda}_i = -\,{{\partial H}\over{\partial x_i}} = -\,\sum_{k=1}^n\,
\lambda_k\,{{\partial f_k}\over{\partial x_i}},\qquad i = \overline{1,n};
\eqno(1.9)$$

\item[(2)] equations of motion (1.1) with respect to the variables $x$ in the form of
Euler equations
$$\dot{x}_i = {{\partial H}\over{\partial\lambda_i}} = f_i\qquad i = \overline{1,n};
\eqno(1.10)$$

\item[(3)]  the universality condition
$$\bigl(\lambda\,\Delta x - H\,\Delta t\bigr)\,\big|_{t_0}^{t_1} = 0.\eqno(1.11)$$
\end{itemize}
\smallskip

It is more believed that the right sides of equations (1.9), (1.10) satisfy the
existence and unique conditions for the Cauchy problem in a limited region $D,$
namely, they satisfy the continuous vector-function $f(x,t)$ and limited elements of the Jacobi matrix $f_x(x,t).$
These equations therefore fit a unique continuous integral trajectory in the
expanded space $R^n \times R$ that issues out of the point with coordinates
$x_0,\lambda_0$ at the initial moment of time~$t_0.$

Note that Euler-Lagrange equations (1.9), \ $\dot{\lambda} = -\,f_x\lambda$,
$f_x = \partial f/\partial x,$  are the Euler equations of
variational Lagrange problem (1.2) with the fixed ends that serve to finding an~extremal of the functional $S,$
$$L_x={d\over dt}L_{\dot{x}},
$$
where $$L_x = {{\partial L}\over{\partial x}} = -\,f_x \lambda,\qquad L_{\dot{x}} =
{{\partial L}\over{\partial\dot{x}}} = \lambda.$$

It is clear that the system of equations (1.9), (1.10) is well-evident to form a~canonical system with
the Hamiltonian $H,$ coordinates $x_i$, and factors $\lambda_i.$ This system is different
from a~familiar canonical system by no more than the linear dependence of $H$ on $\lambda_i$. What is
more, for a total derivative of the function $H = H(x,\lambda,t)$ with respect to the time, we,
according to (1.9) and (1.10), have
\begin{eqnarray*}
{{d H}\over{dt}} &=& {{\partial H}\over{\partial x}}\,\dot{x} + {{\partial
H}\over{\partial\lambda}}\,\dot{\lambda} + {{\partial H}\over{\partial t}} \\
&=&\ {{\partial H}\over{\partial x}}\,{{\partial H}\over{\partial\lambda}} -
{{\partial H}\over{\partial\lambda}}\,{{\partial H}\over{\partial x}} +
{{\partial H}\over{\partial t}} = {{\partial H}\over{\partial t}}.\\
\end{eqnarray*}
    Because by relations (1.5) the expression
$${{dH}\over{dt}} = \dot{\lambda}f + \lambda\dot{f} = \dot{\lambda}f + \lambda f_x\,
\dot{x} + \lambda f_t = \lambda f_t$$
is valid, we shall from here find that on the extremal obeying the necessary conditions for
the extremal the equality $f_x = \partial f/\partial x,\ f_t = \partial f/\partial t$ takes place,
we from here find that, the equality
$$H=h+\int_{t_0}^{t_1}\lambda f_t\, dt, h=const \eqno(1.12) $$ must be
fulfilled.
If the vector function $f$  does not explicitly depend on $t$, it follows
evidently from (1.12) that $H = h.$

We pointed out that from equations (1.9), (1.10) their
contingency follows, where the vector $\lambda(t)$ was that of conjugate variables. But other than
these, the functions $L$  and $H$ are conjugate characteristic functions derived
from the Legendre transformations because, by (1.4), (1.5),
and (1.10) on the extremal trajectories of functional (1.2) while differential
relation (1.1) is provided, the following equalities are resulting:
$$L = {{\partial L}\over{\partial\dot{x}}}\,\dot{x} - H = \lambda\,{{\partial H}
\over{\partial\lambda}} - H = \lambda f - \lambda f = 0.$$

We also further on need certain results of variational analysis concerning the notion
of integral invariance. The calculus of variation [9]
for problem (1.2) with the Lagrangian $L$ the class of curves $C$ with the ends
given on the interval $[t_0,t_1]$ is considered.
It is necessary to find a minimum of the integral $S(C)$ (1.2) along $C$, where $x(t)$ is
the representation of curve $C$. The integral is here taken over
the interval of time $[t_0,t_1]$, within which the curve $x(t)$ is defined.

It is stated that if $M$ is the class of curves with the given ends $A=x(t_0)$, $B=x(t_1)$, and $C_0$
 is a curve from this class, the quantity
$S(C_0)$ is their minimum $S(C)$ with respect to $C\in M$ on the assumption that there exists
a total derivative $\Phi$ of the function $\Psi =\Psi (x, \lambda, t)$
such that $L=\Phi$ along the curve $C_0$ and $L\geq \Phi $ along all other curves from the class $M$.

Explain that we here have $\Phi = \dot{\Psi} = \Psi_x\,\dot{x} + \Psi_\lambda\,\dot{\lambda}
+ \Psi_t.$ The statement in question may be proved straightforwardly: $\forall
\,C \in M$  the value $S(C)$ is as follows: $$S(C) \ge \int_{(C)}\,\Phi\,dt = \int_{(C_0)}\,L\,dt = S(C_0).$$

\begin{definition}
The vector function $g=g(x,t)$ is called a geodesic inclination if there exists a total derivative
$\Phi = \dot{\Psi} = \Psi_x\,\dot{x}
+ \Psi_\lambda\,\dot{\lambda} + \Psi_t$ of the function $\Psi = \Psi(x,\lambda,t)$
such that $\forall\,(x,\dot{x},\lambda,t),$ where $(x,t)$ are lying in the region of definition of the function $g(x,t)$,
the inequality $$L \ge \Phi$$ is valid, and $L = \Phi$ at $\dot{x} = g.$ Here we have $g(x,t):\ R^n\times R \to R^n.$

Thus, $\min\,(L - \Phi) = 0$ is achieved at $\dot{x} = g.$ Providing the coupling equation (1.1), we therefore have
$$L = \lambda\,(\dot{x} - f) = \Psi_x\,\dot{x} + \Psi_\lambda\,\dot{\lambda} +
\Psi_t,$$
and we from where obtain a chain of equalities
$$\Psi_x = L_{\dot{x}},\qquad \Psi_\lambda = L_{\dot{\lambda}} = 0,\qquad \Psi_t
= L - L_{\dot{x}}\,\dot{x} = -\,H.\eqno(1.13)$$
\end{definition}

Let us add to definition 1 [9] that if $g$ is the geodesic inclination, the
geodesic inclination curves are called solutions of the differential
equations $\dot{x}
= g,$ where $g = g(x,t).$ The family of such curves is called a geodesic
family. The geodesic family is named the geodesic flow if it covers the region,
where the function $g$ is only once defined.

Note that in the agreed notation, the integral
$$S = \int_{(C)}\,\bigl[\,L(x,g,\lambda,t) - (\dot{x} - g)\,L_g(x,g,\lambda,t)\,
\bigr]\,dt\eqno(1.14)$$
named the invariant Hilbert integral, where $L_g = \partial L/\partial
g,$ is independent of the integration
path $C,$ but only dependent on the end of the curve $C.$ Due to relations
(1.13) at $\dot{x} = g$, integral (1.14) is indeed equal to the value
$\int_{(C)}\Phi\,dt = \Psi(B) - \Psi(A),$ where $A$ and $B$ are the
beginning end of the integration path; in our case $\int_{(C)}\Phi\,dt = \Psi(B) - \Psi(A).$

With the use of inequality of the form $L - \Phi \ge 0$ and the substitution of $\Phi$ as a
integrand (1.14), we obtain the Weierstrass condition $E\ge 0$, where the
Weierstrass function $E = E(x,g,\dot{x},\lambda,t)$ is of the form
$$E = L(x,\dot{x},\lambda,t) - L(x,g,\lambda,t) - (\dot{x} - g)\,L_g(x,g,\lambda,t).
\eqno(1.15)$$

The $\min\,(L - \Phi) = 0$ is thus achieved at $\dot{x} = g.$
When relation equation (1.1) is satisfied, we therefore have
$$L = \lambda\,(\dot{x} - f) = \Psi_x\,\dot{x} + \Psi_\lambda\,\dot{\lambda} +
\Psi_t,$$   we from where obtain the chain of equalities
$$\Psi_x = L_{\dot{x}},\qquad \Psi_\lambda = L_{\dot{\lambda}} = 0,\qquad \Psi_t
= L - L_{\dot{x}}\,\dot{x} = -\,H.\eqno(1.13)$$

Add to the definition 1 [9] that if $\dot{x}= g$ is the geodesic
inclination, the solutions of the geodesic equation $\dot{x}=0$, where
 $g = g(x,t),$
are referred to as the geodesic inclination curves. The family of such curves is called
a geodesic family. The geodesic family
is referred to as a geodesic flow if it only once covers the region of definition for the function
 $g = g(x,t)$.

 Point out that in the agreed notation the integral
$$S = \int_{(C)}\,\bigl[\,L(x,g,\lambda,t) - (\dot{x} - g)\,L_g(x,g,\lambda,t)\,
\bigr]\,dt\eqno(1.14)$$ named the Hilbert integral, where
$L_g = \partial L/\partial g$, is independent of the integration path  $C,$ but dependent on the
ends of the curve $C.$ By (1.13), the integral (1.13) at  $\dot{x} = g$ is indeed equal to
$\int_{(C)}\Phi\,dt = \Psi(B) - \Psi(A)$ at $\dot{x} = g$, where $A$ and $B$ are the initial and terminal
ends of integral pathes, being $L(x,g,\lambda,t) = \lambda\,(g - f)$ in our case.

Using the inequality $L - \Phi \ge 0$ from definition 1 and substituting the integrand (2.14)  for $\Phi$,
we obtain the Weierstrass condition $E \ge 0,$ where the Weierstrass function
$E = E(x,g,\dot{x},\lambda,t)$ has the form
$$E = L(x,\dot{x},\lambda,t) - L(x,g,\lambda,t) - (\dot{x} - g)\,L_g(x,g,\lambda,t).
\eqno(1.15)$$

Point to the fact that if $g$ is the geodesic inclination, we have $L = \Phi + E,$
hence, if $\forall\,C \in M,$ it joints two levels $\Psi = \Psi_1$ and $\Psi = \Psi_2,$
we from here obtain the Weierstrass formula of calculus of variation
$$S(C) = \Psi_2 - \Psi_1 + \int_{(C)}\,E\,dt,$$
or
$$S(C) = S(C_0) + \int_{(C)}\,E\,dt,$$
where $C_0$ is the curve of flow, and $C$ is any curve joining the levels of the curve
$C_0.$

A necessary Weierstrass condition  $(E \ge 0)$ for the minimum of the action functional $S$
(see (1.2)) can be probed as follows. When $g=g(x,t)$ substitutes for the function $\dot{x}$,
the Lagrange function $L(x,\dot{x},
\lambda,t)$ changes into the function $L(x,g,\lambda,t).$ Expanding into series, we
find $$L(x,\dot{x},\lambda,t) = L(x,g,\lambda,t) + {{\partial L}\over{\partial\dot{x}}}
\,\bigg|_{\dot{x}=g}\,(\dot{x} - g) + \alpha\,\bigl(|\,\dot{x} - g\,|\bigr),$$
where $\alpha(\cdot)$ denotes a set of infinitesimals of the order higher than
first order with respect to $|\,\dot{x} - g\,|.$ The written expression gives
 the Weierstrass function $E$ (1.15) if we neglect~$\alpha.$

    We have
for the system (1.1) in the corresponding variational problem, according to formulas
(1.4) and (1.5), $$L = \lambda\,\bigl[\,\dot{x} - f(x,t)\,\bigr],\qquad {{\partial L}\over{\partial
\dot{x}}} = \lambda.$$ Hence, $$E = \lambda\,(g - f) - \lambda(\dot{x} - f) -
\lambda\,(g - \dot{x}) \equiv 0,$$
i.~e. the Weierstrass condition is satisfied.

Insert into consideration a scalar function $U = U(x,\lambda,t)$
doubly continuously differentiable
with respect to its arguments for the Hamilton system (1.9), (1.10), which we
name the controlling function. That the function $U(x,\lambda,t)$ is supposedly not given
beforehand, but found from quite definite target conditions.

The controlled mapping from the old variables $x(t) \in R^n,$ $\lambda(t)
\in R^n$ to the new ones $y(t) \in R^n,\,\mu(t) \in R^n$ is realized
using the relation
$$x + U_\lambda = y,\qquad \lambda - U_x = \mu,\eqno(1.16)$$
where
$U_\lambda = \partial U/\partial\lambda,\,U_x = \partial U/\partial x.$
Note that the record $x = q,\,\lambda = p, \,y = Q,\,\mu = P$
is usual, further we however follow the established notation.
In addition, it from expressions (1.16) follows that in case of identical
transformations $x = y,\,\lambda = \mu$, we have
$U_\lambda \equiv 0,\,U_x \equiv 0,$ and $U \equiv 0$ may
be taken as the controlling function  $U = U(x,\lambda,t)$.

Give attention to the fact that the mapping $(x,\lambda) \to (y,\mu)$
given by the equalities (1.16) is the generalization of the known
Hilbert-Courant mapping [5] for the parametric generative
function $\Gamma(a,b,t)$ in a canonical system ($a$
 and $b$  are here the vector parameters) at the change from
the variables  to $q = q(a,b,t),\,p = p(a,b,t)$ to those of
$Q = Q(a,b,t),\,P = P(a,b,t)$, by the rule
$$q = a - {1\over 2}\,\Gamma_b,\quad p = b +
{1\over 2}\,\Gamma_a,\quad Q = a +{1\over 2}
\,\Gamma_b,\quad P = b - {1\over 2}\,\Gamma_a,
\eqno(1.17)$$
where $\Gamma_a = \partial\Gamma/\partial a,\,\Gamma_b =
\partial\Gamma/\partial b.$ The removing
parameters $a$ and $b$ give in these formulas clearly give
rise to the relations of form (1.16).

An important advantage of mapping (1.16) over the parametric form of mapping
(1.17) consists in the lack of any intermediate parameters and in the subsequent
necessity for their determination.

    The Jacobian of mappings (1.16) is assumed to be non-zero:
$$\det\,\biggl({{\partial y}\over{\partial x}}\biggr)= {{D(y_1,\,...,
\,y_n)}\over{D(x_1,\,...,\,x_n)}} \not = 0,\nonumber$$
$$\det\,\biggl({{\partial\mu}\over{\partial\lambda}}\biggr)= {{D(\mu_1,\,...,\,
\mu_n)}\over{D(\lambda_1,\,...,\,\lambda_n)}} \not= 0. \eqno(1.18)$$

In relations (1.18) there are determinants of the matrices $E \pm \partial^2
U/\partial x \partial\lambda$ respectively, where
$E$ is a single-valued matrix of $n$ dimensions.

Instead of relation (1.16), a controlled sympletic mapping may be taken that
resembles by its form the Hilbert-Courant mapping
$$x + {1\over 2}\,I U_\lambda = y,\qquad \lambda - {1\over 2}\,I U_x = \mu,\qquad
I = \left(\matrix{0&E\cr -\,E&0\cr}\right),\eqno(1.19)$$
where $E$  is the sympletic matrix. It is easy to see that in the
 presence of
(1.19), the Jacobians are the same and equal to
$$\det\,\biggl({{\partial y}\over{\partial x}}\biggr) = \det\,\biggl({{\partial
\mu}\over{\partial\lambda}}\biggr) = 1 + {1\over 4}\,\det\,\biggl({{\partial^2
U}\over{\partial x\,\partial\lambda}}\biggr).$$
We however confine ourselves in our further arguments to mappings (1.16).
\bigskip
\section{Essential traits of controlled mapping}
The leading question that interests us at the given stage is
what qualitative changes are introduced by the controlled mapping (1.16) into the
description of canonical system (1.9), (1.10)
$$\dot{x}=H_\lambda,\;\; \dot\lambda =-H_x
$$
with the Hamiltonian dependent on the canonical arguments $x,\lambda, t.$

We have known that the invariant Hilbert integral (1.14) calculated over the
extremal transforms into the action functional $S$  that takes a minimum value, see
principle of the least action in the Hamilton form
$$S = \int_{t_0}^{t_1}\,L(x,\dot{x},\lambda,t)\,dt = \int_{(C_0)}\,L\,dt,$$ independent
of the integral path, but dependent only on the ends of this path. We
have over the extremal
$$S = \int_{(C_0)} L\,dt = \int_{(C_0)}\biggl({{\partial L}\over{\partial\dot{x}}}
\,\dot{x} - H\biggr)\,dt = \int_{(C_0)}\lambda\,dx - H\,dt.\eqno(2.1)$$
The integral (2.1) is the above treated Hilbert integral. The differential
expression  $\lambda\,dx - H\,dt$ being its member has the form of an integral
Poincare-Kartan invariant.

We assume that in the expanded phase space  $R^{2n+1}$ with coordinates
 $x \in R^n,\,\lambda \in R^n,\,t \in R$ and the Hamilton
function $H = H(x,\lambda,t)$  constructed for the system (1.9), (1.10), the controlled transformation
(1.16) is diffeomorphic, where the diffeomorphism is considered as a mutually
single-valued and mutually differentiable mapping.

Let us construct 1-form $\omega^1 = \lambda\,dx - H\,dt$ and refer to the following definition [10, p.~859].

\begin{definition}
Diffeomorphism retaining the external differential form $\omega^2 =
 \sum_{i=1}^n\,d\lambda_i \wedge dx_i$ is
referred to as a canonical mapping.
\end{definition}

As the phase trajectories (phase flow) $(x,\lambda)$ of system (1.9), (1.10), representing
[1] the rotor lines for the form $\lambda\,dx - H\,dt$ such that by the Stokes
lemma ($\oint_{\gamma_1}\,\omega^1 = \oint_{\gamma_2}\,\omega^1,$)
the statement (the theorem on the integral Poincare-Kartan invariant)
follows that
$$\oint_{\gamma_1}\,\lambda\,dx - H\,dt = \oint_{\gamma_2}\,\lambda\,dx - H\,dt,$$
where $\gamma_1,\gamma_2$ are the closed curves enveloping one and the same tube of phase (integral)
trajectories of Hamilton system (1.9), (1.10).

In these relations the form $\lambda\,dx$ is referred to as a relative integral Poincare
invariant. For the two-dimensional part of the tube of rotor $\sigma$, the Stokes formula
$$\oint_\gamma\,\lambda\,dx = \int\int_\sigma\,d\lambda \wedge dx$$
is valid, the integral invariance then follows of the 2-form $\omega^2 = d\lambda
\wedge dx$ for the phase flow $(x,\lambda)$~[1].

The canonical transformations convert the Hamilton system into a~Hamiltonian system as well.
We in this connection turn our attention to a~known theorem.

\begin{theorem}\ Let $F:\ R^{2n} \to R^{2n}$  be the canonical transformation of phase space that converts
the point with coordinates $(x,\lambda)$ into that with coordinates $(y,\mu).$ In the new phase
coordinates $(y,\mu)$, canonical equations (1.9), (1.10) then have the canonical form
$${{dy}\over{dt}} = {{\partial G}\over{\partial\mu}},\qquad {{d\mu}\over{dt}} =
-\,{{\partial G}\over{\partial y}}\eqno(2.2)$$
with the former Hamilton function
$$G(y,\mu,t) = H(x,\lambda,t).\eqno(2.3)$$
\end{theorem}

This theorem is proved with the use of a~canonical property of mapping, namely:
$$\oint_\gamma\,\lambda\,dx = \oint_\gamma\,\mu\,dy$$
over any closed curve $\gamma,$ from where
$$\oint_\gamma\,\lambda\,dx - \mu\,dy = 0$$
and, therefore, the integral $\int_{(x_0,\lambda_0)}^{(x_1,\lambda_1)}\,\lambda\,dx - \mu
\,dy = Q,$ where $\int_{(x_0,\lambda_0)}^{(x_1,\lambda_1)}\,\lambda\,dx - \mu
\,dy = Q$ depends on no path of integration, but
an initial $(x_0,\lambda_0)$ and final $(x_1,\lambda_1)$ point of path. Here we have
$$\lambda\,dx - \mu\,dy = dQ.$$
$Q = Q(x,\lambda,t)$  is here an arbitrary continuously differentiable function of its arguments.
This equality in an expanded phase space $R^{2n+1}$  transforms into the equality
$$\lambda\,dx - H\,dt = \mu\,dy - H\,dt + dQ.\eqno(2.4)$$

The trajectories of canonical system (1.9), (1.10) are represented by rotor
lines of the form $\mu\,dy - G\,dt + dQ$ (see details in [1]).
Comparing this form with form (2.4) in the right side, we therefore arrive at
the conclusion that equality (2.3), $G(y,\mu,t) = H(x,\lambda,t)$, is valid.
If the Hamilton function does not change in the canonical transformation,
i.~e. $G = H,$ such a transformation calls quite a canonical transformation.

If we have the considered case of controlled mapping (1.16), $y = y(x,\lambda,
t),\,\mu = \mu(x,\lambda,t),$ using
formula (2.4), where on the right there is  $G$ instead of $H$, we arrive at an
important conclusion that in the canonical transformation of phase space
dependent on time $t$, the canonical equations (1.9), (1.10) in the variables $y, \mu, t$ have
a canonical form (2.2) with a new Hamilton function
$$G(y,\mu,t) = H(x,\lambda,t) + {{\partial Q}\over{\partial t}},\eqno(2.5)$$
where, as before, we have  $(x_0,\lambda_0)$ at the fixed initial point
$$Q(x_1,\lambda_1,t) = \int_{(x_0,\lambda_0)}^{(x_1,\lambda_1)}\,\lambda\,dx -
\mu\,dy$$.

So, for the Poincare-Kartan invariant to retain its form, we must require
that the changed part of this integral invariant is a total differential:
$$\lambda\,dx - H\,dt = \mu\,dy - G\,dt + {{\partial Q}\over{\partial x}}\,dx +
{{\partial Q}\over{\partial\lambda}}\,d\lambda + {{\partial Q}\over{\partial t}}
\,dt,$$
from where we obtain the system of relations whose fulfillment must be required
for providing the canonicity of the controlled transformation:
$$G - H = Q_t,\qquad (\lambda - Q_x)\,dx = \mu\,dy + Q_\lambda\,d\lambda.$$

To use the results of Theorem~1, it needs to demonstrate what conditions
are required for controlled transformation (1.16) to be canonical.

\begin{theorem}
\ \ \ \ Let conditions $(1.18)$ be fulfilled and, in addition, the controlling
function $U = U(x,\lambda,t)$ and the variable $\lambda(t) \in
R^n$ in transformation (1.16) satisfy the differential
equality
$$(U_x - \lambda)\,dU_\lambda = U_\lambda\,d\lambda,\eqno(2.6)$$
where $dU_\lambda = U_{\lambda x}\,dx + U_{\lambda\lambda}\,d\lambda + U_{\lambda t}
\,dt.$ The transformation (1.16) of variables $x,\lambda \to y,\mu$ will then be canonical. This
transformation changes Hamilton system (1.9), (1.10) with the Hamiltonian $H(x,\lambda,t)$ into
the Hamilton system (2.2) with the Hamiltonian $G(y,\mu,t)$ by rule (2.5), where $Q(x,\lambda,t) \equiv U(x,\lambda,t).$
\end{theorem}
{\bf Proof.} The canonicity criterion similar to canonicity criterion
given in [2,~3] actually follows from (2.4), (2.5): if in the transformation of variables $x,\lambda \to y,\mu$ the
differential form
$$W = \mu\,dy - \lambda\,dx - (G - H)\,dt \eqno(2.7)$$
is a total differential for a function  $V = V(x,\lambda,t),$ i.~e. $W =
dV,$ where $dV = V_x\,dx + V_\lambda\,d\lambda + V_t\,dt,$ we have a
canonical transformation.

In the considered case, the arguments $x,\lambda,y,\mu$ in the Hamiltonians $H$ and $G$ are related
through expressions (1.16), and the difference $G - H$ satisfies equality (2.5).
We must in this way take the function $-\,Q$ as a function $V$, following equality (2.4).
Setting further $V = -\,Q = -\,U,$ we obtain the following expression for the form $W$ (2.6):
$$W = d\,(-\,U).\eqno(2.8)$$

Let us verify the fulfillment of this condition. Substitute relations (1.16) and
(2.5) into the form $W$ (2.7). We then get
$$W = (\lambda - U_x)\,d\,(x + U_\lambda) - \lambda\,dx - U_t\,dt = (\lambda -
U_x)\,dU_\lambda - U_x\,dx - U_t\,dt.$$
With regard for the theorem condition (2.6), we obtain for expression (2.8)
$$W = -\,U_\lambda\,d\lambda - U_x\,dx - U_t\,dt = -\,dU.$$
It has thereby been argued that the controlled mapping (1.16) is
a~canonical mapping that establishes the theorem itself.

Some arguments of a general nature, which are immediate from Theorem 2 that
deals with the canonicity of the controlled mapping, can be proved. Note that
the scalar condition (2.6) may be written in a form all the more compact
$$U_x\,dU_\lambda = d\,(\lambda U_\lambda).$$
The condition (2.6) is in any case the equal of the canonicity criterion on
the existence of some function $V(x,\lambda,t)$ such as $W = dV,$ where the form $W$ is
described by equality (2.7).

How could we take advantage of canonicity condition (2.6)? In principle,
criterion (2.6) can be seen as the criterion that is satisfied by the controlling
 function $U(x,\lambda,t)$  (or rather its vector differential components $U_x,U_\lambda$)
 and
vector-function of variables $\lambda(t).$

However, to limit the choice of the function $U,$ on the assumption of condition
(2.6), where the vector $\lambda$ is given by equation (1.9), is a rather strict
requirement. Because it first of all groups together $2n$ unknown
vector-functions $U_x$ and $U_\lambda.$ The condition (2.6) would therefore be appropriate for
the choice of a not uniquely defined vector $\lambda$ rather than in the choice of $U$.

Now demonstrate how this can be done. Let the Jacobi matrix as a preliminary be
denoted as $f_x = A,\,A =
A(x,t).$ Then it is known [11] that the solution of vector differential
equation (1.9) $\lambda(t_0) = \lambda_0$  may be written as
$$\lambda (t)=B(x,t)\lambda_0,
$$ where $B(x,t)$ is the matrix obeying the matrix differential equation
$$\dot{B}(x,t) = -\,A(x,t)\,B(x,t),\qquad B(x_0,t_0) = E,$$
where $x_0 = x(t_0)\, E$ is the unit $n$-matrix. The formula (2.9) for the
mapping $\lambda_0 \to \lambda$
determines the differential homomorphism of class $C^1.$

Substitute expressions (1.9) and (2.9) into canonicity criterion (2.6):
$$(\lambda - U_x)\,dU_\lambda = U_\lambda A\lambda\,dt,$$
from where we arrive at the scalar equation
$$\bigl(\dot{U}_\lambda - A^*U_\lambda\bigr)\,\lambda = U_x\dot{U}_\lambda,$$
or
$$\bigl(\dot{U}_\lambda - A^*U_\lambda\bigr)\,B\lambda_0 = U_x\dot{U}_\lambda,$$
where * over indicates the transposition operation and all vectors are
multiplied as scalars.

Let us also denote the vector $C = B^*\,\bigl(\dot{U}_\lambda - A^*U_\lambda\bigr).
$ Then the latter equation is
$$C\lambda_0 = U_x\dot{U}_\lambda,$$
where $C\lambda_0=(C,\lambda_0)=\sum_{i=1}^n C_i\lambda_{0i}$
may be considered as an equation with respect to one, for example, $k$-th
initial condition $\lambda_{0k}:$
$$\lambda_{0k}={{U_x\dot U_\lambda -\sum_{i-1}^nC_i\lambda_{0i}}\over C_k}\mid_{t=t_0},
$$
where  at the top of the formula there is $i=\overline{1,n},$ $i\neq k$, and $t_0$ is the fixed initial moment of
time.

One more point needed to be made that the controlling function $U$, by the
before made premises and especially the assumption that $U = Q,\ U(x_1,\lambda_1,
t) = \int_{(x_0,\lambda_0)}^{(x_1,\lambda_1)}\,\lambda\,dx - \mu\,dy$, is of the nature
of potential (energetic) function. In this connection the function $U$ may be
named the controlling potential function, or controlling potential.

It is significant that the canonical condition (2.6) is realized
with a feedback. If the initial moment of time $t_0$ to be really considered as a
fixed one, the expression for $\lambda_{0k}$  is determined, having regard to the solutions of
the corresponding equations of dynamics (1.9), (1.10), i.~e. to the
action-functional extremals, through the current values of phase variables. We
arrive by that at the problem of canonical synthesis as far as here we construct the law
for forming the initial value  $\lambda_{0k},$ which provides the fulfillment of the canonicity
condition for a controlled mapping and covers all manifold of other initial
data.

Let us show that the above obtained condition of canonicity and the
Hamilton--Jacobi equation for the controlling function are closely related.

\begin{theorem}\ \ \ \ Let canonicity criterion (2.6) be fulfilled for mapping (1.16). Then,
in order that mapping (1.16) to be a solution of the Cauchy problem for Hamilton
equations (2.2)
$$\dot{y} = {{\partial G}\over{\partial\mu}},\qquad \dot{\mu} = -\,{{\partial G}
\over{\partial y}},\qquad y_0 = y(t_0),\qquad \mu_0 = \mu(t_0),$$
the controlling function $U$ must satisfy the Hamilton--Jacobi equation with the
Hamiltonian~$G:$
$${{\partial U}\over{\partial t}} = G\biggl(x + {{\partial U}\over{\partial\lambda}},
\,\lambda - {{\partial U}\over{\partial x}},\,t\biggr),\eqno(2.10)$$
where  $U = U(x,\lambda,t),\,U(x,\lambda,t_0) = 0.$
\end{theorem}
{\bf Proof.} To prove the theorem, we need equality (2.5) given at $Q = U.$ If
canonicity criterion (2.6) is fulfilled, we have
$${{\partial U}\over{\partial t}} + H(x,\lambda,t) = G(y,\mu,t).\eqno(2.11)$$

It remains to be noted that in the Cauchy problem for the system of canonical
equations (2.2), the variables $x$ and $\lambda$ became initial points
of the canonical system $$x = y(t_0),\qquad \lambda = \mu(t_0)$$
 for the trajectory $y = y(t),\,\mu = \mu(t)$ with the Hamiltonian $G(y,\mu,t)$.  As this
takes place, the motion equations are shaped into the simplest form that
appropriates to a zero Hamiltonian $H:\ H(x,\lambda,t) = 0.$ Therefore,
 to determine the controlling
function $U$ from equation (2.11), we shall obtain equation (2.10). This completes the proof.

It is ready to discover that for the Cauchy problem of system (1.9), (1.10) at
$x_0 = x(t_0),\,\lambda_0
= \lambda(t_0).$ Theorem 3 can be reformulated to terms of an old Hamiltonian
$H(x,\lambda,t)$ for Cauchy's problem of system (1.9), (1.10).
We in this case have $y = x(t_0),\,\mu = \lambda(t_0)$, and $G(y,\mu,t) = 0$ in equation (2.11).

\begin{theorem}
For the Hamilton system (1.9), (1.10), the Cauchy problem solutions
may be presented by canonical mappings (1.16) when the function $U$ satisfies the
Hamilton--Jacobi equation with the Hamiltonian  $H$ having the form
$${{\partial U}\over{\partial t}} + H(x,\lambda,t) = 0,$$
where, as before, $U = U(x,\lambda,t),\,U(x,\lambda,t_0) = 0,$  and the variables
 $x$ and $\lambda$ are solutions of equations (1.9),
(1.10) respectively.
\end{theorem}
{\bf Remarks}\begin{itemize}\item[(1)] We can see from the structure of the proofs that Theorems 3 and 4
are indeed invertible, i.~e. represent the necessary and sufficient conditions
for the solutions of corresponding Hamilton--Jacobi controlling function
equations (2.10), (2.11) to exist.

\item[(2)] Draw our attention to the fact that the functions $U = U(x,\lambda,t)$ are common
in their notation in equations (2.10) and (2.11),  but different in values.

\item[(3)] The controlling function $U$  should not be confused with the generating
function that is actively used in the classical formalism of the theory of
Hamilton equations, canonical transformations, and integration of differential
equations by the Hamilton--Jacobi method.
\end{itemize}

These functions are different in meaning of their formation and further
application. Their main distinction consists in the dependence of the
controlling function on the old variables $x,\lambda$  only, including time $t.$ Recall
that the generating function $\tilde{U}$ is an arbitrary function of mixed (old and
new) variables.

Of course, the function $U$  can be given the nature of generating function $\widetilde{U}.$ With
this aim in view, let us write mappings (1.16) in a general form,
$$y = \varphi(x,\lambda,t),\qquad \mu = \psi(x,\lambda,t),\eqno(2.12)$$
where $\varphi \equiv x + U_\lambda,\ \psi \equiv \lambda - U_x,$ assuming
 that at $\partial^2\varphi/\partial x\,\partial\lambda \not= 0,\,\partial^2\psi
/\partial x\,\partial\lambda \not= 0$ they are solvable in the old variables
$$x = x(y,\mu,t),\qquad \lambda = \lambda(y,\mu,t).\eqno(2.13)$$

The substitution of one of formulas (2.13) into the function $U$  for $x$  or $\lambda$
(there can be nothing but four variants) leads to the appearance of a
generating function $\widetilde{U}.$ The advantage in the use of controlled mapping (1.16)
with the controlling function $U$ over the mapping with the generating function $\widetilde{U}$
is obvious: the use of $U$ does not provide for the resolution (reversibility) of
equations (2.12) in $x$ and $\lambda$ and for the conversion to the explicit dependences (2.13).

It was specified in [1] that the generating function formalism seems to be
'depressive in its non-invariance and essentially uses phase-state coordinate
structure'. Taking into account this note, let us study the question on
invariant properties at controlled mapping (1.16): $(x,\lambda,
t) \to (y,\mu,t).$

Consider, for example, the converted form of the integral Poincare-Kartan
invariant in the canonical conversion $(x,\lambda) \to (y,\mu)$
under the generating function $\widetilde{U}(x,\mu,t)$
$$\lambda\,dx - H\,dt = \mu\,dy - G\,dt + d\widetilde{U},\eqno(2.14)$$
where $d\widetilde{U} = \widetilde{U}_x\,dx + \widetilde{U}_\mu\,d\mu + \widetilde{U}_t
\,dt.$ It follows in particular that
$$G = \widetilde{U}_t + H,\qquad \lambda = \widetilde{U}_x,\qquad y = \widetilde{U}_\mu.$$

The canonical transformation $(x,\lambda) \to (y,\mu)$ obtained with the generating function $\widetilde{U}(x,\mu,t)$ is
admissible if the condition
$$\det\,{{\partial^2\widetilde{U}}\over{\partial x\,\partial\mu}} \not= 0$$
has been fulfilled in relation (2.14). It is obvious that this condition depends on the choice of
new canonical variables.

On the contrary, for the controlling function $U(x,\lambda,t)$ with the converted form
$$\lambda\,dx - H\,dt = \mu\,dy - G\,dt + dU,$$
the existence condition for the canonical transformation directly depends on
the condition $\det\,\bigl(E \pm \partial^2 U/\partial x\,\partial\lambda\bigr) \not=
0$ (1.18) applied to the controlled mappings (1.16). This
condition is invariant with respect to the new canonical variables and
may be affected by the canonical change of variables.

Consider further the Lagrangian $K$ derived from the Lagrangian $L(x,\dot{x},\lambda,
t)$ by means of
controlled change of variables (1.16) if we set in addition that relations
(1.16) are resolvable over the old variables   $x$ and $\lambda$ by formulas (2.13):
$$L(x,\dot{x},\lambda,t) = L\bigl[\,x(y,\mu,t),\,\dot{x}(y,\mu,t),\,\lambda(y,
\mu,t),\,t\,\bigr]\ $$
$$=\ K(y,\dot{y},\mu,\dot{\mu},t),\eqno(2.15)$$
where $\dot{x}(y,\mu,t) = x_y\dot{y} + x_\mu\dot{\mu} + x_t.$
Here we have $$\Omega = x_y = {{\partial x}\over{\partial y}},\qquad x_\mu = {{\partial x}\over
{\partial\mu}},\qquad x_t = {{\partial x}\over{\partial t}}$$ in standard notation.

\begin{definition} The Lagrangian $K$ (2.15) obtained after the change of variables
in the Lagrangian $L,$ will be  called an induced Lagrangian. The corresponding Euler
equation for $K$ will be called the induced Euler equation and its solution
the induced extremals.
\end{definition}

The following theorem on the invariance of mappings (1.16), (2.13) is valid,
which generalizes a known statement of calculus of variation [9].

\begin{theorem}
Let all suppositions about the existence of mappings (1.16), (2.13)
be fulfilled. Then the extremals for the variational problem with the
Lagrangian $L$ are the induced extremals for the variational problem with the induced
Lagrangian $K.$
\end{theorem}
{\bf Proof.} Let us use the notation $v = K_{\dot{y}},$ i.~e.
$$v=\frac{\partial K}{\partial\dot{ y}}=\frac{\partial L}{\partial\dot{x}}
\frac{\partial\dot{x}}
{\partial y}=\lambda\Omega,\eqno(2.16)$$
as far as $\dot{x}_{\dot{y}} = x_y.$ We find now $K_y:$
$${\partial K\over\partial y}={\partial L\over\partial x}{\partial x\over\partial y}
+{\partial L\over\dot{x}}{\partial\dot{x}\over y}
+{\partial L\over\partial\lambda}{\partial\lambda\over\partial y}\eqno(2.17)$$
$$=L_x\Omega+\lambda\Xi +(\dot{x}-f)\Lambda,$$
where the notation
$$\Xi = {{\partial\dot{x}}\over{\partial y}},\qquad \Lambda = {{\partial\lambda}
\over{\partial y}}$$
has been used.

Let us find the value of matrix $\Xi$ from expression (2.17):
$$\Xi = {{\partial\dot{x}}\over{\partial y}} = {\partial\over{\partial y}}\bigl(
x_y\dot{y} + x_\mu\dot{\mu} + x_t\bigr) = {{\partial^2 x}\over{\partial y\,\partial y}}
\,\dot{y} + {{\partial^2 x}\over{\partial y\,\partial\mu}}\,\dot{\mu} +
{{\partial^2 x}\over{\partial y\,\partial t}}.$$
On the other hand, the value of matrix $\dot{\Omega},$ where
$\Omega\ =\ x_y = \Omega(y,\mu,t)$, may be written as
$$\dot{\Omega} = {{\partial\Omega}\over{\partial y}}\dot{y}
+ {{\partial\Omega}\over
{\partial\mu}}\,\dot{\mu} + {{\partial\Omega}\over{\partial t}} =
{{\partial^2 x}\over{\partial y\,\partial y}}\,\dot{y} + {{\partial^2 x}\over
{\partial\mu\,\partial y}}\,\dot{\mu} + {{\partial^2 x}\over{\partial t\,\partial
y}}.$$ It from here follows  that we have the equality $\Xi = \dot{\Omega}.$\\

Then we construct an induced Euler equation with the induced Lagrangian
$K$ using equations (2.16), (2.17):
\begin{eqnarray*}
 {d\over{dt}}\,{{\partial K}\over{\partial\dot{y}}} - {{\partial K}
\over{\partial y}} = \dot{v} - K_y\ \
=\ {d\over{dt}}\,(\lambda\Omega) - L_x\Omega - \lambda\dot{\Omega} + (\dot{x} -
f)\,\Lambda\\
\qquad\qquad\qquad\qquad =\ \bigl(\dot{\lambda} - L_x\bigr)\,\Omega + (\dot{x} - f)\,\Lambda.
\end{eqnarray*}
It thus follows that the equations $\dot{\lambda} = L_x,\ \dot{x} = f,$ and consequently
$\dot{v} = K_y$, take place on the
extremals of the Euler equations generated by the Lagrangian $L$. That establishes
the theorem.

Note in passing that we would arrive at the conclusion which has been formulated in
Theorem~5 that mappings (1.16), (2.13) are invariant if we use the invariance in
the Poincare-Kartan form when the choice of variables to be canonical.

The representation of the controlled mapping is, of course, not be limited by
formulas (1.16). They were entered for definiteness sake and for the demonstration of the
relation with Hilbert-Courant mapping (1.17). The same treatment may be applied
very well to the controlled mapping of the form
$$x \pm U_\lambda = y,\qquad \lambda \pm U_x = \mu,\eqno(2.18)$$
or
$$x \pm U_x = y,\qquad \lambda \pm U_\lambda = \mu\eqno(2.19)$$
may be as readily entered, the signs on the left parts of equations (2.18), (2.19)
can be either identical or opposite. Still more exotic combinations seem to be
possible when constructing new variables $y,\mu.$

The main requirement for the new variables,
nevertheless, as well as for the old ones, is their canonicity. Take, for example, a controlled mapping of the form
$$x + U_x = y,\qquad \lambda - U_\lambda = \mu.\eqno(2.20)$$

Satisfy that the transformations (2.20) lead in specific situations to the canonical
variables $y,\mu$. In this case, a theorem analogous to Theorem~2 may be proved.

\begin{theorem}Let conditions (1.18) be fulfilled with the controlling function $U = U(x,\lambda,t)$
such that its partial derivatives $U_x,U_\lambda$ satisfy equations (2.20). If
the equality
$$(\lambda - U_\lambda)\,dU_x = (U_\lambda - U_x)\,dx - U_\lambda\,d\lambda
\eqno(2.21)$$
is then satisfied,
the transformation of the variables $x,\lambda \to y,\mu$ will be canonical. System (2.2) has here
the Hamiltonian $G(y,\mu,t)$, (2.5), where $Q(x,\lambda,t)
\equiv  U(x,\lambda,t).$
\end{theorem}
{\bf Proof.} Theorem~6 schematically appears in the following form. We formulate
differential form (2.7), where $W=dV,\;\; V=V(x,\lambda, t)$ is any function. Further it needs to check that
if the canonicity criterion (2.21) is fulfilled, the function $V=-U$, where $U = U(x,\lambda,t)$  is the
controlling function, it is a desired one.
We have by mapping (2.20) for form (2.7)
\begin{eqnarray*}W = (\lambda - U_\lambda)d (x + U_x) - \lambda\,dx - U_t\,dt\
\end{eqnarray*}
\begin{eqnarray*}
= \lambda\,dU_x - U_\lambda\,dx - U_\lambda\,dU_x - U_t\,dt
\end{eqnarray*}
\begin{eqnarray*}
= (\lambda - U_\lambda)\,dU_x - U_\lambda\,dx - U_t\,dt
\end{eqnarray*}
\begin{eqnarray*}
= -\,U_x\,dx - U_\lambda\,d\lambda - U_t\,dt = d\,(-\,U).
\end{eqnarray*}
The latter record is, obviously, completes the proof of the theorem.

In addition, the canonicity criterion (2.21) can be provided by the choice of
a certain $k$-th initial condition $\lambda_{0k}.$ In this connection to solve equation (1.9), we
take equality (2.9). Considering initial condition (1.1), criterion (2.21) can
be written in the following form:
$$(\lambda - U_\lambda)\,dU_x = \bigl[\,(U_\lambda - U_x)\,f + U_\lambda A\lambda
\,\bigr]\,dt,$$
hence, we get the equation
$$(\lambda - U_\lambda)\,\dot{U}_x = (U_\lambda - U_x)\,f + U_\lambda A\lambda,$$
or
$$\bigl(\dot{U}_x - A^*U_\lambda\bigr)\,B\lambda_0 = U_\lambda\dot{U}_x +
(U_\lambda - U_x)\,f.$$

An expression for $\lambda_{0k}$ remains to be written. Resolving the latter scalar equation in
$\lambda_{0k},$ we obtain the formula
$$\lambda_{0k} = {{U_\lambda\dot{U}_x + (U_\lambda - U_x)\,f - \sum_{i=1}^n\,
D_i\,\lambda_{0i}}\over{D_k}}\,\bigg|_{t=t_0},$$
where $i \not= k,\ k = \overline{1,n},\ D = B^*\,\bigl(\dot{U}_x - A^*U_\lambda\bigr).$

\noindent{\bf Example~1} The controlled mapping (2.2) was chosen not accidentally. It is convenient with its help
to make rotations of the phase coordinate system $(x,\lambda).$

For example, the rotation of axes  $x$ and $\lambda$ through the right angle in an anticlockwise
direction that is executed by transformations (2.20) must be consistent with
the equalities $(y = \lambda,\,
\mu = -\,x):\ x + U_x = \lambda,\,\lambda - U_\lambda = -\,x,$ from where
the expressions for $U_x$ and $U_\lambda$ follow:
$$U_x = \lambda - x,\qquad U_\lambda = \lambda + x.$$
The controlling function
$$U = U(x,\lambda,t) = {{\lambda^2}\over 2} - {{x^2}\over 2} + \lambda x +
u(t),\qquad \lambda^2 = (\lambda,\lambda),\quad x^2 = (x,x)$$
obviously provides an appropriate rotation of axes $x$ and $\lambda$ for every continuously
differentiable function $u(t)$.

\section{Energy transformations and controlling fields}
By the controlling field is meant a scalar functional field in $R$ that is given
by the values of function $U(x,\lambda,t).$ We previously mentioned that the controlling
function $U = U(x,\lambda,t)$ was of energy nature from the assumptions made on the canonicity
of a~controlled mapping. Turn to this point.

Take once more expression (1.2) for the action functional $S$. The action
functional, where $t=t_1$ in integral (1.2), $t$ being the running time, and the initial point $(x_0,\lambda_0,t_0)$ fixed, is said to be an action function
$$S(x,\lambda,t) = \int_{t_0}^t\,L(x,\dot{x},\lambda,t)\,dt.\eqno(3.1)$$

Try to discover a relation between the action function $S$ (3.1) and the controlling
function $U$. If such a relation will be estimated, we shall be able to form an algorithm
for determining the controlling function.

Denote through $\gamma$ the extremal joining its initial point
$(x_0,t_0)$ with its terminal at
the running point $(x,t)$. For the action function $S(x,\lambda,t)$ (3.1), we thus have as well
$$S(x,\lambda,t) = \int_\gamma\,L(x,\dot{x}\lambda,t)\,dt.$$

According to [1], we may demonstrate that if the initial point $(x_0,t_0)$ is fixed, the
differential of the action function $S(x,\lambda,t)$ has the form of integral Poincare-Kartan
invariant (cf. with expression (2.1)):
$$dS = \lambda\,dx - H\,dt,\eqno(3.2)$$ where
$$dS = {{\partial S}\over{\partial x}}\,dx + {{\partial S}\over{\partial\lambda}}
\,d\lambda + {{\partial S}\over{\partial t}}\,dt,\qquad \lambda = {{\partial L}
\over{\partial\dot{x}}},$$
and the quantity $H = \lambda\dot{x} - L$ is determined in view of a finale
velocity $\dot{x}$ of the trajectory $\gamma.$

The proof of equality (3.2) is based on lifting the extremal $\gamma$ from the space
$(x,t)$ to the expanded phase space $(x,\lambda,t)$, where $\lambda = \partial
L/\partial\dot{x}$. In this case, the extremal is substituted
by a phase trajectory that is among the variety of rotor lines having the form
$\lambda\,dx - H\,dt.$

By the fact that the relations
$${{\partial S}\over{\partial t}} = -\,H(x,\lambda,t),\qquad \biggl({{\partial S}
\over{\partial x}} - \lambda\biggr)\,dx + {{\partial S}\over{\partial\lambda}}\,
d\lambda = 0\eqno(3.3)$$ follow from equality (3.2), we conclude that
the action function $S = S(x,\lambda,t)$ satisfies the
Hamilton-Jacobi equation
$${{\partial S}\over{\partial t}} + H\biggl(x,\,{{\partial S}\over{\partial x}},\,
t\biggr) = 0.\eqno(3.4)$$ The form (1.2) can indeed be written as
$$dS = L\,dt,$$
or $\dot{S} = L,$ where $L = L(x,\dot{x},\lambda,t) = \lambda\,(\dot{x} - f).$

We have then
$$S_x\,\dot{x} + S_\lambda\,\dot{\lambda} + S_t = L,$$
from where the equalities
$$S_x = L_{\dot{x}} = \lambda,\qquad S_\lambda = L_{\dot{\lambda}} = 0,\qquad
S_t = L - L_{\dot{x}}\,\dot{x} = -\,H = -\,\lambda f$$
follow (compare with equalities (1.13)). With an additional condition
$S_\lambda = 0$, the first equation in system (3.3) is therefore converted to form (3.4).

Having regard to the written relations, the Hamilton--Jacobi equation (3.4) of
the scalar action function $S(x,\lambda,t)$  can then be presented in the form
$${{\partial S}\over{\partial t}} = -\,{{\partial S}\over{\partial x}}\,f(x,t).$$

The system of $2n$  equations
$${{\partial S}\over{\partial x}} = \lambda,\qquad {{\partial S}\over{\partial
\lambda}} = 0$$
has in this case to be thought of as a system of equations in $2n$ unknowns $x(t)$ and
$\lambda(t)$.

If the Cauchy problem with the initial condition $S(x,\lambda,
t_0) = S_0(x,\lambda)$ is formulated for equation
(3.4), its solution is then reduced to the solution of canonical Hamilton
equations (1.9), (1.10):
$$\dot{x} = {{\partial H}\over{\partial\lambda}},\qquad \dot{\lambda} = -\,{{\partial
H}\over{\partial x}}$$
with the initial conditions
$$x(t_0) = x_0,\qquad \lambda(t_0) = \lambda_0 = {{\partial S_0}\over{\partial x}}
\,\bigg|_{x=x_0}.\eqno(3.5)$$

The solution of this problem on the space $(x,t)$ comprises, as we know, the extremal
$x = x(t)$ for the given variational principle $\Delta\,\int L\,dt = 0$.
 This extremal issues out of the initial
point $x_0$ and is called the characteristic of problem (3.4). We integrate equality
(3.2) along the characteristic that joins the points $A_0 = (x_0,t_0)$ and $A = (x,t)$. And for the action
function $S$  with the initial condition $S_0$, we find the expression
$$S(A,\lambda) = S_0(A_0,\lambda_0) + \int_{A_0}^A\,L(x,\dot{x},\lambda,t)\,dt,
\eqno(3.6)$$
 which gives the solution of the problem (3.4). Note in addition that the initial
condition $\lambda_0$ must be correlated to the initial condition $\lambda_{0k}$ (3.5).

Let us briefly consider the most important autonomous case. Let time $t$ now not
come explicitly in the expression of Hamilton function $H$, i.~e. $\partial H/\partial t = 0$. The equality
$dH/dt = \partial H/\partial t$  was above obtained in view of he Hamilton equations (1.9), (1.10). We
therefore have $dH/dt = 0$ in the given variant, from where it follows that  $H = h = {\rm const}$  is the first integral
(Jacobi integral) of equations (1.9), (1.10).

By the Legendre transformation, $H = \bigl(\partial L/\partial\dot{x}\bigr)\,\dot{x}
- L = \lambda f$ and $f = f(x)$. The Legendre function $L$ is then independent
of time too; $L = L(x,\dot{x},\lambda) = \lambda\,(\dot{x} - f)$ and $\partial H/\partial
t = -\,\partial L/\partial t = 0$.

Let the surface $H(x,\lambda) = h$ be projected from the expanded phase space $(x,\lambda,t)$  to the space $(x,\lambda)$. In
this case, the time $t$ actually does not vary ($(dt = 0)$), the total variation $\Delta$ changes into
the isochronous variation $\delta$, and
the expression (3.2) takes
a shortened form as compared to the relative integral
Poincare invariant:
$$d\,S=\lambda\,d\,x.
\eqno(3.7)$$ Hence, the trajectories of Hamilton system (1.9),(1.10) are extremals for the
variational principle corresponding to form (3.7).

The phase trajectories of canonical equations (1.9),(1.10) lying on
the surface $H(x,\lambda) = h$  are in this way extremals of the integral
$$S = \int_\gamma\,\lambda\,dx,\eqno(3.8)$$ and these extremals joint the points $x_0$ and $x_1$.
The formulated principle forms the
contents of the Maupertuis--Lagrange least (stationary) action principle
validated by Lagrange [12] with regard to the Lagrange action  $J$, namely,
$$J = \int_{t_1}^{t_2}\,2T\,dt,$$
where $T = mv^2/2$ is the kinetic energy of a material point of mass $m$
and velocity $v = ds/dt$, namely
$$J = \int_{t_1}^{t_2}\,mv^2\,dt = \int_{t_1}^{t_2}\,mv\,{{ds}\over{dt}}\,dt =
\int_{s_1}^{s_2}\,mv\,ds.$$

It is clear that in the agreed notation the Lagrange action $J$ will be equal to
the shortened action integral $S$ (4.8) if the impulse $mv = \lambda$ and distance $s = x$ are taken.

Let the controlled mapping (1.16) be in the autonomous case given by the
controlling function $U = U(x,\lambda)$. If this mapping is quite a canonical one, then by Theorem 1,
the conversion from Hamilton equations (1.9),(1.10) to those of form (2.2) takes
place according to equality (2.3) for the Hamiltonians $G = H,$ where $G = G(y,\mu),\,H = H(x,\lambda)$.

In so doing, the relation $\lambda\,dx - \mu\,dy = dQ$ is valid (see expression (2.4), where $Q = Q(x,\lambda))$. We have
by that
$$\lambda\,dx = \mu\,dy + {{\partial Q}\over{\partial x}}\,dx + {{\partial Q}
\over{\partial\lambda}}\,d\lambda,$$
we get from where
$$(\lambda - Q_x)\,dx = \mu\,dy + Q_\lambda\,d\lambda.$$
Select as a function $Q:\ Q = U,\,U = U(x,\lambda)$.  Then for the canonicity of
transformations (1.16), the equality
$$(\lambda - U_x)\,dx = (\lambda - U_x)\,d\,(x + U_\lambda) + U_\lambda\,d\lambda$$
needs  to be required, or
$$d\,(\lambda U_\lambda) = U_x\,d U_\lambda,$$
i.~e. we again arrive at the canonicity criterion (2.6).

We shall return to basic relation (2.4), where $Q
(x,\lambda,t) \equiv U(x,\lambda,t)$ and the function $U(x,\lambda,t)$ specifies
controlled mappings (1.16). These mappings are canonical when criterion (2.6) is
fulfilled, i.~e. at a certain choice of the initial vector $\lambda_0$.

Further, let us use the concepts of Lagrangian, Hamiltonian, action function, canonical
variables, and write energy the equality (2.4) being main in the controlled mapping
method in the form
$$L\,dt = K\,dt + dU,\eqno(3.9)$$ where $L\,dt = \lambda\,dx - H\,dt,\;L = L(x,\dot{x}\lambda,t)$
 is an initial Lagrange function, $x, \lambda$ are initial canonical variables,
$H = H(x,\lambda,t)$ is the initial Hamilton function,
$K\,dt = \mu\,dy - G\,dt,\,K = K(y,\dot{y}\mu,t)$ is a new Lagrange function:
$$K = \mu\,(\dot{y} - g),\qquad \mu = K_{\dot{y}},\qquad G = \mu g.\eqno(3.10)$$
$G = G(y,\mu,t) = \mu g$ is here a new Hamilton function, where $g = g(y,t)$,
which is connected with the function $K$
by a Legendre transformation and $y,\mu$ are new canonical variables in the new
$(2n + 1)$-dimensional expanded phase space $(y,\mu,t)$.

The function $K(y,\dot{y},\mu,t)$ is corresponded by a new action function
$R = R(y,\mu,t)$:
 $$R(y,\mu,t) = \int_{t_0}^t\,K(y,\dot{y},\mu,t)\,dt = \int_\delta\,K(y,\dot{y},
\mu,t)\,dt,\eqno(3.11)$$ where $\delta$ is the extremal of new variational principle
$\Delta\,\int K\,dt = 0$, $\delta$ being an induced extremal.
 The new canonical equations
with a new Hamiltonian $G$ have the form of equations (2.2)
$$\dot{y} = {{\partial G}\over{\partial\mu}} = g,\qquad \dot{\mu} = -\,{{\partial
G}\over{\partial y}} = -\,g_y\,\mu,\eqno(3.12)$$
where $g = g(y,t)$  is a certain function continuously differentiable in its
arguments.

It is important to note that the choice of the new Hamiltonian $G = G(y,\mu,t)$
and of the
right parts in canonical equations (3.12), respectively (for example for a more
convenient integration or reduction of the systems), in a more simple form (its
normalization) is a main purpose condition of the formation of the controlling
function  $U = U(x,\lambda,t)$.

Using the formalism of action functions, we write relation (3.9) that determines
the controlling potential $U(x,\lambda,t)$ in the form
$$dS = dR + dU,\qquad S - R = U,\eqno(3.13)$$
from where we derive the expression for $U(x,\lambda,t)$:
$$U(x,\lambda,t) = \int_\gamma\,\lambda\,dx - H\,dt - \biggl(\int_\delta\,\mu\,
dy - G\,dt\biggr)\  \eqno(3.14)$$
$$=\ \int_{(x_0,\lambda_0)}^{(x,\lambda)}\,\lambda\,dx - H\,dt - \biggl(\int_{
(y_0,\mu_0)}^{(y,\mu)}\,\mu\,dy - G\,dt\biggr) = \int_{t_0}^t\,\bigl(\lambda\dot{x}
- \mu\dot{y} + G - H)\,dt.$$

In formula (3.14) by the canonicity of controlled transformations, the
written line integrals do not lye
 on the path of integration but depends only
on the finale values $(x,\lambda)$, $(y,\mu)$ at the fixed initial values $(x_0,
\lambda_0),\,(y_0,\mu_0)$.

In addition, in (3.14) the new canonical variables $y,\mu$ are connected with the old
canonical variables $x,\lambda$ by controlled mappings (1.16). Note also that for the new
action function $R$ (3.11) (from equation (3.13)) there takes place a
Hamilton-Jacobi equation of the form
$$R_t = -\,G,\qquad R_y = \mu,\qquad R_\mu = 0,$$
where $G = G(y,R_y,t).$

\noindent {\bf Example~2} Consider an autonomous dynamic system that is square integrable and
described by the equation
$$\dot{x} = f(x),\qquad x \in R^n$$
with the given initial condition $x_0 = x(t_0)$ and suitable requirements to the
vector-function $f(x)$. It is now necessary, using controlled transformations of the
form (1.16), to reduce it to the given system presented by the equation
$$\dot{y} = g(y),\qquad y \in R^n,$$
where the vector-function $g(y)$ is determined from the given Hamilton function $G(y,\mu)$.

Use all previous notation and constructions to solve this problem. We have
$$L(x,\dot{x},\lambda) = \lambda\,(\dot{x} - f),\qquad H(x,\lambda) = \lambda f
= h,$$
where $h > 0$  is a constant equal to the value of total system energy. Canonical
equations in the variables $x,\lambda$ have the form
$$\dot{x} = {{\partial H}\over{\partial\lambda}} = f(x),\qquad \dot{\lambda} =
-\,{{\partial H}\over{\partial x}} = -\,f_x(x)\,\lambda.$$
It is considered that the initial data $x_0,\lambda_0$  are given and $\lambda_0$ is coordinated with
the canonicity criterion (2.6).

Let the new Hamiltonian $G$ meet the requirement
$$G(y,\mu) = a\mu,$$
where $a \in R^n$  is a given constant vector. Due to the autonomy, $H
= G = h$. The new canonical
equations (3.12) may be written in the form
$$\dot{y} = {{\partial G}\over{\partial\mu}} = a,\qquad \dot{\mu} = -\,{{\partial
G}\over{\partial y}} = 0,$$
from where it follows that  $y(t) = at + b,\,\mu(t) = c$.
Here $a,b,c$ are the vectors of known components. Thus,
we have $g(y) = a,\,G(y,\mu) = ac = h.$

It is required, starting from these purpose conditions, to determine the
controlling function $U = U(x,\lambda)$ such that it is reduced to these new Hamilton
equations. Let there be controlling mappings (1.16) and by the choice of $\lambda_0)$
canonicity criterion (2.6) be provided.

Here with regard to (3.9)--(3.14), we have for the situation of finding the
controlling function $U(x,\lambda):$
$$\lambda\,dx = \mu\,dy + dU,\qquad U = S - R,$$
from where, since $H = G = h$, we get a chain of equivalent relations
$$U(x,\lambda) = \int_\gamma\,L\,dt - \int_\delta\,K\,dt\ $$
$$=\ \int_\gamma\,\lambda\,\bigl[\,\dot{x} - f(x)\,\bigr]\,dt - \int_\delta\,c\,
(\dot{y} - a)\,dt\ $$
$$=\ \int_{(x_0,\lambda_0)}^{(x,\lambda)}\,\lambda\,dx - \int_{(y_0,c)}^{(y,c)}
\,c\,dy = \int_{(x_0,\lambda_0)}^{(x,\lambda)}\,\lambda\,dx - c\,(y - y_0),$$
where $y = y(t),\,y_0 = y(t_0)$,  and by the canonicity of transformations, all written line integrals do
not depend on the path of integration, but depend only on the finale values at the
fixed initial data.

Substituting the expression for $y$  (1.16) into the latter equality, we
obtain an equation for the required function $U(x,\lambda)$:
$$U(x,\lambda) = \int_{(x_0,\lambda_0)}^{(x,\lambda)}\,\lambda\,dx - c\,\bigl(x
+ U_\lambda\bigr) + cy_0.$$

This equation is a first-order partial linear differential equation with respect
to the unknown function $U(x,\lambda).$  This has the form
$$U(x,\lambda) + c\,{{\partial U(x,\lambda)}\over{\partial\lambda}} = F(x,\lambda),
\eqno(3.15)$$
where the notation is $$c\,{{\partial U(x,\lambda)}\over{\partial\lambda}} = \sum_{i=1}^n\,c_i\,{{\partial
U(x,\lambda)}\over{\partial\lambda_i}},\quad F(x,\lambda) = \int_{(x_0,\lambda_0)}^
{(x,\lambda)}\lambda\,dx + c\,(y_0 - x),$$
$F(x,\lambda)$ being a known function of variables $x$ and $\lambda.$

An attempt to get an analytical solution of equation (3.15) meets certain
difficulties [13, 14] caused by 2n independent variables $x_1,...,\,x_n,$ $\lambda_1,
...,\,\lambda_n$ and $n$  partial
derivatives $U_{\lambda_1},...,\,U_{\lambda_n}.$ To obtain the relation between $x_1,...,\,x_n,\,\lambda_1,...,\,\lambda_n$
and $U$, we must find a solution
for the corresponding system of ordinary differential equations (characteristic
equations of a specific form) that satisfy the given initial conditions. If we
disregard the theoretical questions concerning the solution of equation (3.15),
we can point to the fact that in the given case this equation can be considered
as a final algorithmic equation for determining the controlling potential that
gives a mapping in the initial phase space $(x,\lambda)$ with subsequent numerical resolution.

Let us apply the obtained results of the method of dynamical controlled-mapping
method to the most important class of mechanical systems, Hamiltonian
systems, which in standard notation of generalized coordinates $q(t) \in R^n$ and generalized
momenta $p(t) \in R^n,$ have the form
$$\dot{q}_i = {{\partial H}\over{\partial p_i}},\qquad \dot{p}_i = -\,{{\partial
H}\over{\partial q_i}},\qquad i = \overline{1,n},\eqno(3.16)$$
where $H = H(q,p,t) = p\dot{q} - L(q,\dot{q},t)$ is the Hamilton function, $L = L(q,
\dot{q},t)$ is the Lagrange function, $L = T - \Pi,\ T$ is the kinetic
energy of the system, $\Pi = \Pi(q,t)$ is its potential energy,  $T = T_2 + T_1 + T_0,\ T_2 = (1/2)
\,\sum\,a_{ij}\,\dot{q}_i\dot{q}_j$ is the quadratic form of
generalized velocities, $A = (a_{ij})$ is the kinetic energy matrix, $T_1 = \sum\,b_i\dot{q}_i$
is a linear form of
generalized velocities, and $T_0 = T_0(q,t)$ is its null-form. Thus, we have
$$L = {1\over 2}\,\sum_{i,j=1}^n\,a_{ij}\,\dot{q}_i\dot{q}_j + \sum_{i=1}^n\,b_i
\,\dot{q}_i + T_0 - \Pi,$$
where the coefficients  $a_{ij},b_i$ depend on the generalized coordinates and time. Hence,
the generalized momenta are expressed by
$$p_i = {{\partial L}\over{\partial\dot{q}_i}} = \sum_{j=1}^n\,a_{ij}\,\dot{q}_j
+ b_i.$$

Since the kinetic energy matrix $A$ is non-degenerated, solving then the latter
linear equation with respect to $\dot{q}$, we get
$$\dot{q}_i = \sum_{j=1}^n\,\alpha_{ij}\,(p_j - b_j),$$
where $A^{-1} = (\alpha_{ij})$. If to substitute this expression into the expression for $H$, we shall
find
$$H = \sum_{i=1}^n\,p_i\dot{q}_i - L = \sum_{i,j=1}^n\,\alpha_{ij}\,(p_i - b_i)
\,p_j\ $$
$$-\ {1\over 2}\,\sum_{i,j,k,l=1}^n\,a_{ij}\,\alpha_{ik}\,\alpha_{jl}\,(p_k - b_k)
(p_l - b_l)$$
$$-\ \sum_{i,j=1}^n\,b_i\,\alpha_{ij}\,(p_j - b_j) - T_0 + \Pi\ $$
$$={1\over 2}\sum_{i,j=1}^n\alpha_{ij}(p_ip_j + b_ib_j) - \sum_{i,j=1}^n
\alpha_{ij}b_ip_j - T_0 + \Pi.$$

In case when the Lagrange function does not depend on time $(\partial L/\partial t = 0)$, we have the
first integral (generalized energy integral, Painleve--Jacobi
integral) as follows:
$$\sum_{i=1}^n\,\dot{q}_i\,{{\partial L}\over{\partial\dot{q}_i}} - L = {\rm const},$$
or $T_2 - T_0 + \Pi = {\rm const}$. If the system is conservative, $T_1 = T_0
= 0$  and the generalized energy potential is
the same as the total energy $T_2\ +\ \Pi
= {\rm const}$ (the total energy conservation law in the
conservative systems).

The Hamilton equations (3.16) are equivalent to the Lagrange equations
$${d\over{dt}}\,{{\partial L}\over{\partial\dot{q}_i}} - {{\partial L}\over{\partial
q_i}} = 0,\qquad i = \overline{1,n},\eqno(3.17)$$
where $\partial L/\partial\dot{q}_i = p_i$ are the generalized momenta, $\partial L/
\partial q_i$ are the generalized forces, the
trajectories of motion for mechanical systems (3.16), (3.17) agreeing with the
extremals of the action functional $\int_{t_0}^{t_1}
\,L(q,\dot{q},t)\,dt$ (principle of least action in the Hamilton
form).

Further, we must rewrite the dynamic relations of mechanical systems in terms
of customary standard notation $x(t)$ and $\lambda(t)$. Let us take for this aim $x(t)
= (q(t),\,\dot{q}(t))$ or $x(t) = (\dot{q}(t),\,q(t))$ and attach the
meaning of vector Hamilton multiplier to $\lambda(t)$. We then resolve equation (3.17)
with respect to a higher derivative and introduce the vector of new derivatives
 $x(t)$. Now we arrive at a normal formulation of the mechanical system in the form of
equation of motion (1.1). The neat process of controlled canonical
transformation of system is completely corresponded to the above described
scheme.

It should be particularly emphasized that the controlled mapping method is
based on the representation the system of equations of motion and Euler
equations as a normal system of differential Cauchy equations (1.9), (1.10)
respectively, giving the Hamilton function $H$ in a special form (1.5), $H = \lambda f.$
Such a formulation sometimes leads to problems.

\noindent{\bf Example~3} In the problem dealing with the determination of Lagrange multipliers
for a ballistic flight (it is a section of path, where an engine is excluded)
in a central gravitation field [7], it is necessary to present the equations of
motion and Euler equations (Euler--Hamilton ones) in a normal Cauchy form and
rewrite the Hamiltonian in the form of (1.5).

The trajectory of motion, as we know from theoretical mechanics, is here represented
by a conic section that in polar coordinates $r,\varphi$ satisfies the equation
$$r = {p\over{1 + e\cos\,(\varphi - \omega)}},$$
where we have the notation: $p$ is a focal parameter ($(p = r$ at\ $\varphi = \pi/2\ +\ \omega),
\,e$ is the
eccentricity, $\omega$  is a polar angle periphery center, i.~e. the trajectory point
closest to the focus of conic section, the angle $a\ =\ \varphi - \omega$ being said to be a true anomaly.
The change of polar angle is determined by the area integral
$$r^2\,\dot{\varphi} = \sigma\,\sqrt{p},$$
where $\sigma = \sqrt{\gamma M},\,\gamma$  is the gravitation constant, $M$ is the mass of a central body.

In order to compile the Euler equations for the multipliers $\lambda,$ we must first
write the equations of motion in the normal form. Let us write them in polar
coordinates. Use some known kinematic expressions for the projections of point
acceleration on the generalized axes of curvilinear polar coordinate system
$$w_r = \ddot{r} - r\dot{\varphi}^2,\qquad w_\varphi = {1\over r}\,{{d\,(r^2
\dot{\varphi})}\over{dt}}.$$

Since the point of mass  $m$ is under the action of only Newton force of
attraction to the fixed center $\vec{F} = -\,m\sigma^2
\vec{r}/r^3$, we obtain the equations of motion
$$\ddot{r} - r\dot{\varphi}^2 = -\,\biggl({\sigma\over r}\biggr)^2,\qquad
{{d\,(r^2\dot{\varphi})}\over{dt}} = 0.$$
Note simultaneously that the second equation leads to an area integral, and the
first one in passing to a new variable $1/r$ and new argument $\varphi$ instead of $t$ leads to
the already known equation of conic sections.

Of course, the equation that describes the motion of a point in a gravitation
field in polar coordinate can also be found with the second-order Lagrange
equations. For this aim, we must write the kinetic equation $T$  and the
generalized forces $Q_r$ and $Q_\varphi$:
$$T = {1\over 2}m\,\bigl(\dot{r}^2 + r^2\dot{\varphi}^2\bigr),\qquad Q_r = -\,m
\,\biggl({\sigma\over r}\biggr)^2,\qquad Q_\varphi = 0$$
and then construct the equations
$${d\over{dt}}\,{{\partial T}\over{\partial\dot{q}_i}} - {{\partial T}\over{\partial
q_i}} = Q_{q_i},\qquad i = 1,2,\quad q_1 = r,\quad q_2 = \varphi$$
being coincident with the above cited equations of motion.

Construct now the equations of motion in a form solvable with respect to derivatives,
i.~e. found in their normal form. To this end, take
the projections of the velocity of a point on the polar axes, namely
$v_r = \dot{r},\,v_\varphi
= r\dot{\varphi}$.  as velocity phase coordinates.
As a result, we obtain the system of four equations of motion
$$\dot{v}_r = {{v_\varphi^2}\over r} - \biggl({\sigma\over r}\biggr)^2,\qquad
\dot{v}_\varphi = -\,{{v_r\,v_\varphi}\over r},\qquad\dot{r} = v_r,
\qquad \dot{\varphi} = {{v_\varphi}\over r}.$$
The Hamilton function $H$   (1.5) then takes the form
$$H = \lambda_1\,\biggl[\,{{v_\varphi^2}\over r} - \biggl({\sigma\over r}\biggr)^2
\,\biggr] + \lambda_2\,\biggl({{-\,v_r\,v_\varphi}\over r}\biggr) + \lambda_3
v_r + \lambda_4\,{{v_\varphi}\over r},$$
and the system of Euler equations for the determination the multipliers $\lambda(t)$  will
be written in the form
$$\dot{\lambda}_1 = \lambda_2\,{{v_\varphi}\over r} - \lambda_3,\qquad \dot{\lambda}_2
= -\,2\lambda_1\,{{v_\varphi}\over r} + \lambda_2\,{{v_r}\over r} - \lambda_4\,
{1\over r},$$
$$\dot{\lambda}_3 = \lambda_1\,{{v_\varphi^2}\over{r^2}} - 2\lambda_1\,{{\sigma^2}
\over{r^3}} - \lambda_2\,{{v_r\,v_\varphi}\over{r^2}} + \lambda_4\,{{v_\varphi}
\over{r^2}},\qquad \dot{\lambda}_4 = 0.$$
 The latter formula can be represented in the form
$\lambda_4 = {\rm const}$. The latter three form  a system of linear equations in
$\lambda_1,\lambda_2,\lambda_3$ with alterable coefficients.
\section{Canonical controlled mappings}
The controlled mappings  (1.16) are constructed on the base of additive occurrence of
controlling functions. The efficiency of their use may change from task to task. Point out in
this connection that  the transformation
 by formulas  (1.16) can be of a more general nature being independent of the choice of initial data.

Along with controlled mapping (1.16)
$$y=x+U_ \lambda , \mu =\lambda-U_x, U=U(x,\lambda),
$$
let us consider a quasi-canonical transformation such that its
vector variables $x$ and $\lambda$ change by an infinitesimal value.
Such a transformation is called the infinitesimal canonical transformation [1, 15].
Suppose
$$y=x+\eta (x, \lambda)\epsilon, \mu=\lambda +\varsigma(x,\lambda)\epsilon, \eqno(4.1)$$
where $\eta (x,\lambda)$,  $\varsigma(x, \lambda)$ are the vector $n$-dimensional functions
of $x,\lambda$, and $\epsilon$ is a small parameter.
If $\eta $ and $\varsigma$ are given arbitrarily, the transformation (4.1) is clearly non-canonical.

We therefore set forth the requirements imposed on the vector-functions
$\eta$ and $\varsigma$ that provide the canonicity
of transformation (4.1), and find a formal analogy between the
vector-functions $U_\lambda, \eta$ and $U_x,\varsigma$.

Let us fall back on equality (2.4) of the form
$$\lambda\, dx=\mu\, dy+dQ,$$
 where $Q=Q(x,\lambda)$ is an arbitrary function of $x$ and $\lambda$.
 Substitute (4.1) into its right side:
$$\lambda\,dx=(\lambda +\epsilon\varsigma)(dx+\epsilon\,d\eta)+dQ.$$

Dropping the terms that contain $\epsilon^2$, we obtain
$$\epsilon (\varsigma\, dx+ \lambda\, \eta)=-dQ,$$
or
$$(\varsigma\, dx -\eta\, d\lambda)+d(\eta\lambda)=-{1\over\epsilon}\, dQ.$$

The latter equality may be written in the form
$$\varsigma\, dx-\eta\, d\lambda =- d\Omega(x,\lambda),$$
where the notation
$$\Omega(x,\lambda)={1\over\epsilon} Q+\lambda\eta$$
is entered for the scalar function $\Omega(x,\lambda)$.

With regard to the fact that $d\Omega =\Omega_x\, dx+\Omega_\lambda\, d\lambda,$
we obtain from (4.2)
$$\varsigma =-\Omega_x,\;\;\; \eta =\Omega_\lambda$$
due to the independence of the quantities $dx$ and $d\lambda$.
Hence, canonical transformation (4.1), being infinitesimal, will be seen as
$$y = x + \Omega_\lambda \varsigma,\qquad \mu = \lambda - \Omega_x \varsigma,
\eqno(4.3)$$
where $\Omega = \Omega(x,\lambda)$  is an arbitrary scalar function of variables $x,\lambda$, which plays a great
role (compared with (1.16)) of a controlling function with a small parameter or,
by another formalism [16--18], of a small controlled perturbation (test
signal).

As for a general non-stationary case, all canonical transformations according
to S.~Lie, see [15], might be defined with one differential condition. Namely, the
reversible transformation $x,\lambda \to y,\mu$ (in this case the functional determinant
$\partial\,(y,\mu)/\partial\,(x,\lambda) \not= 0)$ is
canonical if we have the differential form
$$\lambda\,dx = \mu\,dy + H_0\,dt + dQ,$$
where $H_0$ and $Q$  are arbitrary functions of $4n + 1$ variables
$x,\lambda,y,\mu,t,$ are satisfied identically.

This form should be given as another form that is well known if to subtract the
quantity $H\,dt$  from both parts of equation. Then we shall have
$$\lambda\,dx - H\,dt = \mu\,dy - G\,dt + dQ,\qquad G = H - H_0.$$
It should be noted that independent in these forms will be only $4n + 1$ variables $x,\lambda, y, \mu, t$ from
$2n + 1$ because there exist $2n$
 ratios  between the variables of the form $y = y(x,\lambda,t),\,\mu = \mu(x,\lambda,t)$.

If we apply the canonical controlled transformation (4.3), where $U_\lambda = U_\lambda
(x,\lambda,\varsigma) = \Omega_\lambda(x,\lambda)\,\varsigma$, $U_x = U_x(x,
\lambda,\varsigma) = \Omega_x(x,\lambda)\,\varsigma,$ to a scalar
function  $\psi(x,\lambda)$, then it is possible to find its variation correct to the
magnitudes of the order of $\varsigma^2$ smallness:
$$\psi(x_i,\lambda_i) = \psi\bigl(y_i - U_{\lambda_i},\,\mu_i + U_{x_i}\bigr) =
\psi\bigl(y_i - \Omega_{\lambda_i}\varsigma,\,\mu_i + \Omega_{x_i}\varsigma
\bigr)\ $$
$$=\ \psi(y_i,\mu_i) + \varsigma\,\sum_{i=1}^n\,\biggl(-\,{{\partial\psi}\over
{\partial x_i}}{{\partial\Omega}\over{\partial\lambda_i}} + {{\partial\psi}\over
{\partial\lambda_i}}{{\partial\Omega}\over{\partial x_i}}\biggr),$$
or in the notation of Poisson bracket
$$\bigl\{\,\psi,\,\Omega\,\bigr\} = \sum_{i=1}^n\,\biggl({{\partial\psi}\over
{\partial x_i}}{{\partial\Omega}\over{\partial\lambda_i}} - {{\partial\psi}\over
{\partial\lambda_i}}{{\partial\Omega}\over{\partial x_i}}\biggr),$$
we obtain the formula
$$\psi(x_i,\lambda_i) = \psi(y_i,\mu_i) - \bigl\{\,\psi,\,\Omega\,\bigr\}\,\varsigma
= \psi(y_i,\mu_i) - \bigl\{\,\psi,\,U\,\bigr\}.$$

Interest is provoked by the case of small canonical controlled transformation
when we take the Hamilton function $H(x,\lambda)$ as a function $\Omega$
and infinitesimal interval of time $dt$ as an arbitrary small  parameter.
Suppose we have
$$U(x,\lambda,\varsigma) = \Omega(x,\lambda)\,\varsigma = H(x,\lambda)\,dt.$$

For the infinitesimal canonical transformation (4.3), we may then write
$$y = x + H_\lambda\,dt,\qquad \mu = \lambda - H_x\,dt$$
from where it follows by canonical Hamilton equations (1.9), (1.10) that
$$dx = H_\lambda\,dt,\qquad d\lambda = -\,H_x\,dt$$
according to Hamilton canonical equations (1.9), (1.10)
$$dx=H_\lambda\, dt,\;\; d\lambda=-H_x\, dt.
$$ These relations
actually mean that in a time $dt$ the system $\bigl(
x(t),\,\lambda(t)\bigr)$ will pass to a new state $\bigl(x(t + dt),\,\lambda
(t + dt)\bigr)$, and this
change will be induced by the infinitesimal controlled canonical transformation
such that the basic controlling function of the form
$$\Omega (x,\lambda)=\frac{\partial}{\partial\varepsilon}U(x,\lambda,\varepsilon)
$$
is the Hamilton function $H(x,\lambda)$.  We here assume that the initial function $U$
linearly depends on the parameter $\varsigma$.

In this way, it can be concluded, starting from this analysis, that the movement
of this system is a sequence of infinitesimal controlled mappings, where the
basic controlling function (at a small parameter $\varsigma)$ must be the Hamilton
function itself. Such a result only confirms an important conclusion [19]
that the system can be led to the given specific state by means of sequential
small probe signals (controlled perturbations).

Turn to basic controlled transformations (1.16). We said in the preface to
Section~4 that it was possible to construct canonical controlled mappings
without recourse to specially selected initial data, as was done before.

\begin{theorem}
Suppose the conditions (1.18) are fulfilled for transformation (1.16).
Transformations (1.16) of the variables $x,\lambda \to y,\mu$ are canonical if simultaneously:
\begin{description}
\item{(i)} to choice the vector $U_\lambda$  from the equation
$$\dot{U}_\lambda = A^*\,U_\lambda,\eqno(4.4)$$
\noindent where $A = A(x,t) = f_x(x,t)$;
\item{(ii)} to choice the vector $U_x$ using the orthogonality of vectors $U_x$ and $\dot{U}_
\lambda.$
\end{description}
\end{theorem}
{\bf Proof.} We shall, as before, start from a key equality (2.6) to which controlled
mapping (1.16) leads. If the equality is fulfilled identically,
the transformation (1.16) is canonical.

Using initial Hamilton equation (1.9), we succeed in writing equality (2.6) in
the form
$$\bigl(\dot{U}_\lambda - A^*U_\lambda\bigr)\,\lambda = U_x\dot{U}_\lambda.
\eqno(4.5)$$
The relation (4.5) will be fulfilled identically if we choose the vector $U_\lambda$ such that
equation (4.4) is valid, from where it automatically follows the requirement
that the equality
$$U_x\dot{U}_\lambda = 0\eqno(4.6)$$
to be fulfilled.

The vector differential equation (4.5) is identical in its structure to (1.9) $\dot\lambda=-A\lambda$.
The solution of equation (4.5) is given by the formula
$$U_\lambda = D(x,t)\,U_{\lambda 0},\qquad U_{\lambda 0} = U_\lambda\,\big|_{t=t_0},
\eqno(4.7)$$
where $D(x,t)$ is the $n\times n$-matrix that satisfies the matrix differential equation
$$\dot{D}(x,t) = A^*(x,t)\,D(x,t).\qquad D(x_0,t_0) = E.$$
The theorem is completed.

In particular, it follows from (4.7) that the vector-function $U_\lambda$ is a function
of variables  $x$ and $t$. This means that the controlling function $U(x,\lambda, t)$
is linearly dependent on the $\lambda$-vector elements:
$$U(x,\lambda,t) = \int_{\lambda_0}^\lambda\,D(x,t)\,U_{\lambda 0}\,d\lambda +
u(x,t),$$
where $u(x,t)$ is an arbitrary scalar function of  $x$, and $t$ that can be chosen, for
example, from the orthogonality condition (4.6).

\section*{Conclusions}
The idea to reduce a dynamic system to an expanded Hamilton form with
simultaneous attachment of the formalism of classical calculus of variation was
realized before as well (see, for example, [6--8] and
[20, 21] for adaptive
systems).

It turned to be that such structures were able to work in the construction of
controlled canonical mappings in Hamilton systems for the purpose of writing
the equations of motion in a more comfortable form.

For the purpose of constructing canonical coordinate
mappings in Hamilton systems, the controlling-function method was
developed in the given article, whose leading property was the dependence on the old canonical
variables.

This property is, in the end, allowed additional possibilities in forming an algorithm
for a canonical change of variables, which render it
redundant and less well defined as compared with the algorithm based on the
generating-function method.

Note also a broad range of new variational results, which, being of a general theoretic importance, were arranged as
some theorems. These results may have
practical use in the solution of various problems of canonicity and
study of dynamic Hamilton--Jacobi relations.
\section*{Notation}
 $R^n$ is $n$-dimensional Euclidian space.\\
$\bigtriangleup x_i$, $\delta x_i$ are total and isochronous variations of the component $x_i$ respectively.\\
$\bigtriangleup S$ is a total variation of action functional $S$.\\
$A^*$ is the transpose with respect to the matrix $A$.\\
$S(C)=\int_{(C)}L(x, \dot{x}, \lambda , t)\, dt$ is the integral $S$ along the curve $C$,
where $x=x(t)$ is the representation of curve $C, t\in [t_0,t_1]$.\\
$\{\psi ,\Omega\}$ is the Poisson bracket with the elements $\psi$ and $\Omega$.\\
$B_x=\partial B/\partial x$ is the vector (matrix) of partial derivative of the scalar function (vector) $B$
with respect to the elements of vector $x$.
\newpage
All items are published in Russian.

\end{document}